# Two Finite Classes of Orthogonal Functions

**Mohammad Masjed-Jamei**[a]      **Wolfram Koepf**[b]

[a] *Department of Mathematics, K.N.Toosi University of Technology, P.O.Box 16315-1618, Tehran, Iran,*
E-mail: *mmjamei@kntu.ac.ir , mmjamei@yahoo.com ,*

[b] *Institute of Mathematics, University of Kassel, Heinrich-Plett-Str. 40, D-34132 Kassel, Germany,*
E-mail: *koepf@mathematik.uni-kassel.de,*

**Abstract.** By using Fourier transforms of two symmetric sequences of finite orthogonal polynomials, we introduce two new classes of finite orthogonal functions and obtain their orthogonality relations via Parseval's identity.

**1. Introduction.** There exist ten sequences of real polynomials [6, 10, 11] that are orthogonal with respect to the Pearson distributions family

$$W\left(\begin{array}{cc} d, & e \\ a, & b, & c \end{array}\bigg| x\right) = \exp\left(\int \frac{dx+e}{ax^2+bx+c} dx\right) \quad (a,b,c,d,e \in \mathbf{R}), \qquad (1)$$

and its symmetric analogue [10]

$$W^*\left(\begin{array}{cc} r, & s \\ p, & q \end{array}\bigg| x\right) = \exp\left(\int \frac{r x^2 + s}{x(px^2+q)} dx\right) \quad (p,q,r,s \in \mathbf{R}). \qquad (2)$$

Five of them are infinitely orthogonal with respect to special cases of the two above-mentioned distribution (weight) functions and five other ones are finitely orthogonal [6,10] which are limited to some parametric constraints. The following table shows the main properties of these ten sequences.

*This work has been supported by a grant from "Iran National Science Foundation".*

Table 1: *Characteristics of ten sequences of orthogonal polynomials*

| Symbol | Weight Function | Kind, Interval & Parameters constraint |
|---|---|---|
| $P_n^{(u,v)}(x)$ | $W\begin{pmatrix} -u-v, & -u+v \\ -1, & 0, & 1 \end{pmatrix} x = (1-x)^u(1+x)^v$ | Infinite, $[-1,1]$, $\forall n,\ u>-1, v>-1$ |
| $L_n^{(u)}(x)$ | $W\begin{pmatrix} -1, & u \\ 0, & 1, & 0 \end{pmatrix} x = x^u \exp(-x)$ | Infinite, $[0,\infty)$, $\forall n,\ u>-1$ |
| $H_n(x)$ | $W\begin{pmatrix} -2, & 0 \\ 0, & 0, & 1 \end{pmatrix} x = \exp(-x^2)$ | Infinite, $(-\infty,\infty)$ |
| $J_n^{(u,v)}(x;a,b,c,d)$ | $((ax+b)^2 + (cx+d)^2)^{-u}$ $\times \exp(v \arctan \frac{ax+b}{cx+d})$ | Finite, $(-\infty,\infty)$, $\max n < u - 1/2$, $\det(a,b,c,d) \neq 0$ |
| $M_n^{(u,v)}(x)$ | $W\begin{pmatrix} -u, & v \\ 1, & 1, & 0 \end{pmatrix} x = x^v(x+1)^{-(u+v)}$ | Finite, $[0,\infty)$, $\max n < (u-1)/2$, $v>-1$ |
| $N_n^{(u)}(x)$ | $W\begin{pmatrix} -u, & 1 \\ 1, & 0, & 0 \end{pmatrix} x = x^{-u}\exp(-1/x)$ | Finite, $[0,\infty)$, $\max n < (u-1)/2$ |
| $S_n\begin{pmatrix} -2u-2v-2, & 2u \\ -1, & 1 \end{pmatrix} x$ | $W^*\begin{pmatrix} -2u-2v, & 2u \\ -1, & 1 \end{pmatrix} x = x^{2u}(1-x^2)^v$ | Infinite, $[-1,1]$, $u>-1/2,\ v>-1$ |
| $S_n\begin{pmatrix} -2, & 2u \\ 0, & 1 \end{pmatrix} x$ | $W^*\begin{pmatrix} -2, & 2u \\ 0, & 1 \end{pmatrix} x = x^{2u}\exp(-x^2)$ | Infinite, $(-\infty,\infty)$, $u>-1/2$ |
| $S_n\begin{pmatrix} -2u-2v+2, & -2u \\ 1, & 1 \end{pmatrix} x$ | $W^*\begin{pmatrix} -2u-2v, & -2u \\ 1, & 1 \end{pmatrix} x = x^{-2u}(1+x^2)^{-v}$ | Finite, $(-\infty,\infty)$, $\max n < u+v-1/2$, $u<1/2,\ v>0$ |
| $S_n\begin{pmatrix} -2u+2, & 2 \\ 1, & 0 \end{pmatrix} x$ | $W^*\begin{pmatrix} -2a, & 2 \\ 1, & 0 \end{pmatrix} x = x^{-2a}\exp(-1/x^2)$ | Finite, $(-\infty,\infty)$, $\max n < u-1/2$ |

where the sequence

$$\Phi_n(x) = S_n\begin{pmatrix} r, & s \\ p, & q \end{pmatrix} x = \sum_{k=0}^{[n/2]} \binom{[n/2]}{k} \left( \prod_{i=0}^{[n/2]-(k+1)} \frac{(2i+(-1)^{n+1}+2[n/2])p+r}{(2i+(-1)^{n+1}+2)q+s} \right) x^{n-2k}, \quad (3)$$

is a basic class of symmetric orthogonal polynomials [10] satisfying the equation

$$x^2(px^2+q)\Phi_n''(x) + x(rx^2+s)\Phi_n'(x) - \left(n(r+(n-1)p)x^2 + (1-(-1)^n)s/2\right)\Phi_n(x) = 0.$$
(4)

Except the two finite orthogonal polynomial sequences $S_n\left(\begin{array}{cc}-2u-2v+2, & -2u \\ 1, & 1\end{array}\bigg| x\right)$ and $S_n\left(\begin{array}{cc}-2u+2, & 2 \\ 1, & 0\end{array}\bigg| x\right)$ in table 1, Fourier transforms of all ten aforesaid sequences have been found. Indeed, in [8] Fourier transforms of the generalized ultraspherical polynomials $S_n\left(\begin{array}{cc}-2u-2v-2, & 2u \\ -1, & 1\end{array}\bigg| x\right)$ and generalized Hermite polynomials $S_n\left(\begin{array}{cc}-2, & 2u \\ 0, & 1\end{array}\bigg| x\right)$ have been derived. In [9] the Fourier transform of Routh-Romanovski polynomials $J_n^{(u,v)}(x;a,b,c,d)$ has been obtained, and in [5] the Fourier transforms of finite orthogonal polynomials $M_n^{(u,v)}(x)$ and $N_n^{(u)}(x)$ have been calculated. In this sense, note that the Fourier transforms of classical Jacobi, Laguerre and Hermite polynomials are already known, see e.g. [3,7]. Hence, to complete the analysis of the families of orthogonal polynomials of table 1, only Fourier transforms of the two above-mentioned finite sequences remain, which should be determined. To reach this purpose, we need the general properties of these two sequences.

### 1.1. Finite orthogonal polynomials with weight $x^{-2a}(1+x^2)^{-b}$ on $(-\infty, \infty)$

If $(p,q,r,s) = (1,1,-2a-2b+2,-2a)$ is substituted in (4), then the equation

$$x^2(x^2+1)\Phi_n''(x) - 2x((a+b-1)x^2+a)\Phi_n'(x) + \left(n(2a+2b-(n+1))x^2 + (1-(-1)^n)a\right)\Phi_n(x) = 0,$$
(5)

has the explicit solution

$$\Phi_n(x) = S_n\left(\begin{array}{cc}-2a-2b+2, & -2a \\ 1, & 1\end{array}\bigg| x\right) = \sum_{k=0}^{[n/2]}\binom{[n/2]}{k}\left(\prod_{i=0}^{[\frac{n}{2}]-(k+1)}\frac{2i+2[n/2]+(-1)^{n+1}+2-2a-2b}{2i+(-1)^{n+1}+2-2a}\right)x^{n-2k},$$
(6)

whose monic form is equivalent to the hypergeometric representation

$$A_n^{(a,b)}(x) = \bar{S}_n\left(\begin{array}{cc}-2a-2b+2, & -2a \\ 1, & 1\end{array}\bigg| x\right) = x^n\,_2F_1\left(\begin{array}{c}-[n/2], \ a+1/2-[(n+1)/2] \\ a+b-n+1/2\end{array}\bigg| -\frac{1}{x^2}\right).$$
(7)

$_2F_1(.)$ in (7) is a special case of the generalized hypergeometric functions [1,4]

$$_pF_q\left(\begin{array}{cccc}a_1 & a_2 & \ldots & a_p \\ b_1 & b_2 & \ldots & b_q\end{array}\bigg| x\right) = \sum_{k=0}^{\infty}\frac{(a_1)_k (a_2)_k \ldots (a_p)_k}{(b_1)_k (b_2)_k \ldots (b_q)_k}\frac{x^k}{k!},$$
(8)

where $(r)_k = r(r+1)\ldots(r+k-1)$.
The monic polynomials (7) satisfy the orthogonality relation

$$\int_{-\infty}^{\infty} \frac{x^{-2a}}{(1+x^2)^b} A_n^{(a,b)}(x) A_m^{(a,b)}(x)\, dx = (-1)^n \prod_{j=1}^{n} C_j \begin{pmatrix} -2a-2b+2, & -2a \\ 1, & 1 \end{pmatrix}$$
$$\times \frac{\Gamma(b+a-1/2)\Gamma(-a+1/2)}{\Gamma(b)} \delta_{n,m}, \quad (9)$$

in which

$$C_j \begin{pmatrix} -2a-2b+2, & -2a \\ 1, & 1 \end{pmatrix} = \frac{(j-(1-(-1)^j)a)(j-(1-(-1)^j)a-2b)}{(2j-2a-2b+1)(2j-2a-2b-1)}, \quad (10)$$

and

$$\int_{-\infty}^{\infty} \frac{x^{-2a}}{(1+x^2)^b}\, dx = \int_{0}^{\infty} \frac{t^{-a-1/2}}{(1+t)^b}\, dt = B(-a+\frac{1}{2}; b+a-\frac{1}{2}) = \frac{\Gamma(b+a-1/2)\Gamma(-a+1/2)}{\Gamma(b)}. \quad (11)$$

According to [10], relation (9) is valid only if $m, n = 0, 1, ..., N \le a+b-1/2$ where $N = \max\{m,n\}$; $a < 1/2$; $(-1)^{2a} = 1$ and $b > 0$. Moreover, $B(\lambda_1, \lambda_2)$ in (11) denotes the Beta integral [1] having various definitions as

$$B(\lambda_1; \lambda_2) = \int_0^1 x^{\lambda_1-1}(1-x)^{\lambda_2-1}\, dx = \int_{-1}^{1} x^{2\lambda_1-1}(1-x^2)^{\lambda_2-1}\, dx = \int_0^{\infty} \frac{x^{\lambda_1-1}}{(1+x)^{\lambda_1+\lambda_2}}\, dx$$
$$= 2\int_0^{\pi/2} \sin^{(2\lambda_1-1)} x \cos^{(2\lambda_2-1)} x\, dx = \frac{\Gamma(\lambda_1)\Gamma(\lambda_2)}{\Gamma(\lambda_1+\lambda_2)} = B(\lambda_2; \lambda_1), \quad (12)$$

in which

$$\Gamma(z) = \int_0^{\infty} x^{z-1} e^{-x}\, dx \qquad \text{Re}(z) > 0, \quad (13)$$

denotes the well-known Gamma function satisfying the equation $\Gamma(z+1) = z\Gamma(z)$.

## 1.2. Finite orthogonal polynomials with weight $x^{-2a} e^{-1/x^2}$ on $(-\infty, \infty)$

Similarly, if $(p,q,r,s) = (1, 0, -2a+2, 2)$ is substituted in (4), then the equation

$$x^4 \Phi_n''(x) + 2x((1-a)x^2 + 1)\Phi_n'(x) - \left(n(n+1-2a)x^2 + 1 - (-1)^n\right)\Phi_n(x) = 0, \quad (14)$$

has the explicit solution

$$\Phi_n(x) = S_n \begin{pmatrix} -2a+2 & 2 \\ 1 & 0 \end{pmatrix} x = \sum_{k=0}^{[n/2]} \binom{[n/2]}{k} \left( \prod_{i=0}^{[\frac{n}{2}]-(k+1)} \frac{2i + 2[n/2] + (-1)^{n+1} + 2 - 2a}{2} \right) x^{n-2k},$$
$$(15)$$

whose monic form is equivalent to the hypergeometric form

$$B_n^{(a)}(x) = \bar{S}_n\begin{pmatrix} -2a+2 & 2 \\ 1 & 0 \end{pmatrix} x = x^n {}_1F_1\begin{pmatrix} -[n/2] \\ a+(-1)^n/2 \end{pmatrix} \frac{1}{x^2}. \qquad (16)$$

Moreover, the orthogonality relation corresponding to these polynomials takes the form

$$\int_{-\infty}^{\infty} x^{-2a} e^{-\frac{1}{x^2}} B_n^{(a)}(x) B_m^{(a)}(x) dx = \left( (-1)^n \prod_{j=1}^{n} C_j\begin{pmatrix} -2a+2 & 2 \\ 1 & 0 \end{pmatrix} \right) \Gamma(a-\frac{1}{2}) \delta_{n,m}, \qquad (17)$$

where [10]

$$C_j\begin{pmatrix} -2a+2 & 2 \\ 1 & 0 \end{pmatrix} = \frac{-2(-1)^j(j-a)-2a}{(2j-2a+1)(2j-2a-1)}, \qquad (18)$$

and $m, n = 0, 1, \ldots, N \leq a - 1/2$ with $N = \max\{m, n\}$ and $(-1)^{2a} = 1$.

It is known that some orthogonal polynomials are mapped onto each other by the Fourier transform. In this paper, we follow this approach for the two finite orthogonal polynomials ( described in sections 1.1 and 1.2 ) to obtain two new classes of finite orthogonal functions via Parseval's identity.

## 2. Fourier transform of monic polynomials $A_n^{(a,b)}(x)$ and $B_n^{(a)}(x)$ and their orthogonality relations

The Fourier transform of a function, say $g(x)$, is defined by [2]

$$\mathbf{F}(s) = \mathbf{F}(g(x)) = \int_{-\infty}^{\infty} e^{-isx} g(x) dx, \qquad (19)$$

and for the inverse transform we have

$$g(x) = \frac{1}{2\pi} \int_{-\infty}^{\infty} e^{isx} \mathbf{F}(s) ds. \qquad (20)$$

For $g, h \in L^2(\mathbb{R})$, the Parseval identity related to a Fourier transform is given by [2]

$$\int_{-\infty}^{\infty} g(x) \overline{h(x)} dx = \frac{1}{2\pi} \int_{-\infty}^{\infty} \mathbf{F}(g(x)) \overline{\mathbf{F}(h(x))} ds. \qquad (21)$$

By noting the relations (9) and (21), let us define the functions

$$\begin{cases} g(x) = x^{-2\alpha}(1+x^2)^{-\beta} A_n^{(c,d)}(x) & \text{s.t. } (-1)^{2\alpha} = 1, \\ h(x) = x^{-2l}(1+x^2)^{-u} A_m^{(v,w)}(x) & \text{s.t. } (-1)^{2l} = 1, \end{cases} \qquad (22)$$

in terms of the monic polynomials (7) to which we shall apply the Fourier transform.

Notice that for both above functions the Fourier transform exists. For instance, for $g(x)$ defined in (22) we have

$$\mathbf{F}(g(x)) = \int_{-\infty}^{\infty} e^{-isx}(1+x^2)^{-\beta} x^{-2\alpha} A_n^{(c,d)}(x)\,dx$$

$$= \int_{-\infty}^{\infty} e^{-isx}(1+x^2)^{-\beta} x^{-2\alpha+n} \left(\sum_{k=0}^{[n/2]} \frac{(-[n/2])_k (c+1/2-[(n+1)/2])_k}{(c+d-n+1/2)_k k!} \frac{(-1)^k}{x^{2k}}\right) dx \qquad (23)$$

$$= \sum_{k=0}^{[n/2]} \frac{(-[n/2])_k (c+1/2-[(n+1)/2])_k (-1)^k}{(c+d-n+1/2)_k k!} \left(\int_{-\infty}^{\infty} e^{-isx}(1+x^2)^{-\beta} x^{-2\alpha+n-2k}\,dx\right).$$

Now it remains in (23) to evaluate the definite integral

$$I_{n,k}(s;\alpha,\beta) = \int_{-\infty}^{\infty} e^{-isx}(1+x^2)^{-\beta} x^{-2\alpha+n-2k}\,dx. \qquad (24)$$

There are two ways to compute the integral (24). In the first way, by noting that $(-1)^{2\alpha} = 1$, we can directly compute $I_{n,k}(s;\alpha,\beta)$ for $n = 2m$ as follows

$$I_{2m,k}(s;\alpha,\beta) = \int_{-\infty}^{\infty} \left(\sum_{j=0}^{\infty} \frac{(-isx)^j}{j!}\right) x^{-2\alpha+2m-2k}(1+x^2)^{-\beta}\,dx = \sum_{j=0}^{\infty} \frac{(-1)^j i^j s^j}{j!}\left(\int_{-\infty}^{\infty} x^{j-2\alpha+2m-2k}(1+x^2)^{-\beta}\,dx\right) =$$

$$\sum_{r=0}^{\infty} \frac{(-1)^r s^{2r}}{(2r)!}\left(2\int_0^{\infty} x^{2r-2\alpha+2m-2k}(1+x^2)^{-\beta}\,dx\right) = \sum_{r=0}^{\infty} \frac{(-1)^r s^{2r}}{(2r)!} B(r-\alpha+m-k+\tfrac{1}{2};\beta-r+\alpha-m+k-\tfrac{1}{2}),$$

$$(25)$$

where we have used the third kind of beta integral in (12).
The last sum in (25) can be represented in terms of a hypergeometric function, so

$$I_{2m,k}(s;\alpha,\beta) = \frac{\Gamma(-\alpha+m-k+1/2)\Gamma(\beta+\alpha-m+k-1/2)}{\Gamma(\beta)} {}_1F_2\left(\begin{array}{c} -\alpha+m-k+1/2 \\ 1/2,\ -\beta-\alpha+m-k+3/2 \end{array} \middle| \frac{s^2}{4}\right).$$

$$(26)$$

By knowing that

$$\int_{-\infty}^{\infty} x^{j-2\alpha+2m+1-2k}(1+x^2)^{-\beta}\,dx = 0 \text{ for any } j = 0,2,4,\ldots,$$

this method can be similarly applied to $I_{2m+1,k}(s;\alpha,\beta)$ so that after some computations we obtain

$$I_{2m+1,k}(s;\alpha,\beta) = (-is)\frac{\Gamma(-\alpha+m-k+3/2)\Gamma(\beta+\alpha-m+k-3/2)}{\Gamma(\beta)} {}_1F_2\left(\begin{array}{c} -\alpha+m-k+3/2 \\ 3/2,\ -\beta-\alpha+m-k+5/2 \end{array} \middle| \frac{s^2}{4}\right).$$

$$(27)$$

Hence, combining both relations (26) and (27) and using the identity

$$[\tfrac{n+1}{2}] - [\tfrac{n}{2}] = \frac{1-(-1)^n}{2}, \qquad (28)$$

gives the final form of (24) as

$$I_{n,k}(s;\alpha,\beta) = \Gamma(-\alpha-k+\frac{1}{2}+[\frac{n+1}{2}])\Gamma(\beta+\alpha+k-\frac{1}{2}-[\frac{n+1}{2}])$$
$$\times \frac{(-is)^{\frac{1-(-1)^n}{2}}}{\Gamma(\beta)} {}_1F_2\left(1-\frac{(-1)^n}{2},\ \begin{matrix}1/2-\alpha-k+[(n+1)/2]\\ -\beta-\alpha-k+3/2+[(n+1)/2]\end{matrix}\middle|\frac{s^2}{4}\right). \quad (29)$$

The second way of computing $I_{n,k}(s;\alpha,\beta)$ is that we respectively suppose $n=2m$ and $n=2m+1$ and then directly apply the cosine and sine Fourier transforms [2] to the function $(1+x^2)^{-\beta}x^{-2\alpha+n-2k}$. For example, by noting that $(-1)^{2\alpha}=1$ we have

$$I_{2m,k}(s;\alpha,\beta) = \int_{-\infty}^{\infty} \cos(sx)(1+x^2)^{-\beta}x^{-2\alpha+2m-2k}dx - i\int_{-\infty}^{\infty}\sin(sx)(1+x^2)^{-\beta}x^{-2\alpha+2m-2k}dx$$

$$= 2\int_0^{\infty}\cos(sx)(1+x^2)^{-\beta}x^{-2\alpha+2m-2k}dx$$

$$= \frac{\Gamma(-\alpha+m-k+1/2)\Gamma(\beta+\alpha-m+k-1/2)}{\Gamma(\beta)} {}_1F_2\left(\begin{matrix}-\alpha+m-k+1/2\\ 1/2,\ -\beta-\alpha+m-k+3/2\end{matrix}\middle|\frac{s^2}{4}\right). \quad (30)$$

**Remark 1.** To prove the last equality of (30), one can use dominated convergence theorem (DCT) [12] so that define the sequence

$$\Phi_n(x) = (1+x^2)^{-\beta}x^{-2\lambda}\sum_{k=0}^{n}\frac{(-1)^k(sx)^{2k}}{(2k)!} \quad \text{for}\quad (-1)^{\lambda}=1, \quad (31)$$

and then conclude

$$|\Phi_n(x)| \leq \cosh(sx)(1+x^2)^{-\beta}x^{-2\lambda} = \Phi(x). \quad (32)$$

Inequality (32) allows us to explicitly compute $I_{2m,k}(s;\alpha,\beta)$ in (30) as the same form as we have done in (25).
By referring to remark 1, we can similarly conclude that

$$I_{2m+1,k}(s;\alpha,\beta) = \int_{-\infty}^{\infty}\cos(sx)(1+x^2)^{-\beta}x^{-2\alpha+2m+1-2k}dx - i\int_{-\infty}^{\infty}\sin(sx)(1+x^2)^{-\beta}x^{-2\alpha+2m+1-2k}dx$$

$$= (-2i)\int_0^{\infty}\sin(sx)(1+x^2)^{-\beta}x^{-2\alpha+2m+1-2k}dx$$

$$= (-is)\Gamma(-\alpha+m-k+\frac{3}{2})\Gamma(\beta+\alpha-m+k-\frac{3}{2}) {}_1F_2\left(\begin{matrix}-\alpha+m-k+3/2\\ 3/2,\ -\beta-\alpha+m-k+5/2\end{matrix}\middle|\frac{s^2}{4}\right). \quad (33)$$

Therefore, the result (29) would simplify (23) as

$$\mathbf{F}(g(x)) = \frac{1}{\Gamma(\beta)} \Gamma(-\alpha + \frac{1}{2} + [\frac{n+1}{2}]) \Gamma(\beta + \alpha - \frac{1}{2} - [\frac{n+1}{2}])(-is)^{\frac{1-(-1)^n}{2}} \times$$

$$\sum_{k=0}^{[n/2]} \frac{(-[n/2])_k (c+1/2-[(n+1)/2])_k (\beta+\alpha-1/2-[(n+1)/2])_k}{(c+d-n+1/2)_k (1/2+\alpha-[(n+1)/2])_k k!} \quad (34)$$

$$\times {}_1F_2 \left( \begin{array}{c} 1/2-\alpha-k+[(n+1)/2] \\ 1-\frac{(-1)^n}{2}, \ -\beta-\alpha-k+3/2+[(n+1)/2] \end{array} \bigg| \frac{s^2}{4} \right).$$

If for simplicity in (34) we define the symmetric function

$$A_n(x; p_1, p_2, p_3, p_4) = x^{\frac{1-(-1)^n}{2}} \sum_{k=0}^{[n/2]} \frac{(-[n/2])_k (p_3+1/2-[(n+1)/2])_k (p_1+p_2-1/2-[(n+1)/2])_k}{(p_3+p_4-n+1/2)_k (1/2+p_1-[(n+1)/2])_k k!}$$

$$\times {}_1F_2 \left( \begin{array}{c} 1/2-p_1-k+[(n+1)/2] \\ 1-\frac{(-1)^n}{2}, \ -p_1-p_2-k+3/2+[(n+1)/2] \end{array} \bigg| \frac{x^2}{4} \right),$$

(35)

then clearly

$$\mathbf{F}(g(x)) = (-i)^{\frac{1-(-1)^n}{2}} \frac{1}{\Gamma(\beta)} \Gamma(-\alpha + \frac{1}{2} + [\frac{n+1}{2}]) \Gamma(\beta + \alpha - \frac{1}{2} - [\frac{n+1}{2}]) A_n(s; \alpha, \beta, c, d). \quad (36)$$

By substituting (36) in the Parseval identity (21) and noting (22) one gets

$$2\pi \int_{-\infty}^{\infty} x^{-2(\alpha+l)} (1+x^2)^{-(\beta+u)} A_n^{(c,d)}(x) A_m^{(v,w)}(x) \, dx = i^{\frac{(-1)^n-(-1)^m}{2}} \Gamma(-\alpha+\frac{1}{2}+[\frac{n+1}{2}])$$

$$\times \frac{1}{\Gamma(\beta) \ \Gamma(u)} \Gamma(\beta+\alpha-\frac{1}{2}-[\frac{n+1}{2}]) \Gamma(-l+\frac{1}{2}+[\frac{m+1}{2}]) \Gamma(u+l-\frac{1}{2}-[\frac{m+1}{2}]) \quad (37)$$

$$\times \int_{-\infty}^{\infty} A_n(s;\alpha,\beta,c,d) A_m(s;l,u,v,w) \, ds.$$

Now, if in the left hand side of (37)

$$c = v = \alpha + l \quad \text{and} \quad d = w = \beta + u, \quad (38)$$

then according to the orthogonality relation (9) we have,

**Theorem 1.** *The special function $A_n(x; p_1, p_2, p_3, p_4)$ defined in (35) satisfies the finite orthogonality relation*

$$\frac{1}{2\pi}\int_{-\infty}^{\infty} A_n(x;\alpha,\beta,p,q) A_m(x;p-\alpha,q-\beta,p,q)\,dx = \prod_{j=1}^{n} \frac{\left(-j+(1-(-1)^j)p\right)\left(j-(1-(-1)^j)p-2q\right)}{(2j-2p-2q+1)(2j-2p-2q-1)} \times$$

$$\frac{\Gamma(\beta)\Gamma(q-\beta)\Gamma(p+q-1/2)\Gamma(-p+1/2)\,\delta_{n,m}}{\Gamma(q)\Gamma\left(-\alpha+\frac{1}{2}+[\frac{n+1}{2}]\right)\Gamma\left(\alpha+\beta-\frac{1}{2}-[\frac{n+1}{2}]\right)\Gamma\left(\alpha-p+\frac{1}{2}+[\frac{n+1}{2}]\right)\Gamma\left(p+q-\alpha-\beta-\frac{1}{2}-[\frac{n+1}{2}]\right)},$$
(39)

*for* $m,n = 0,1,\ldots,N = \max\{m,n\} \leq p+q-1/2$, $p < 1/2$, $(-1)^{2p} = 1$, $q > \beta > 0$, $0 < \alpha < 1/2$ *and* $\alpha + \beta > 1/2$.

This approach can similarly be applied for the monic polynomials $B_n^{(a)}(x)$ in (16). First, let us define the functions

$$u(x) = x^{-2a} e^{\frac{-1}{2x^2}} B_n^{(b)}(x) \quad \text{and} \quad v(x) = x^{-2c} e^{\frac{-1}{2x^2}} B_m^{(d)}(x) \quad \text{for} \quad (-1)^{2a} = (-1)^{2c} = 1. \quad (40)$$

The Fourier transform of e.g. $u(x)$ is computed as

$$\mathbf{F}(u(x)) = \int_{-\infty}^{\infty} e^{-isx} x^{-2a} e^{\frac{-1}{2x^2}} B_n^{(b)}(x)\,dx = \int_{-\infty}^{\infty} e^{-isx} e^{\frac{-1}{2x^2}} x^{-2a+n} \left(\sum_{k=0}^{[n/2]} \frac{(-[n/2])_k}{(b+(-1)^n/2)_k} \frac{x^{-2k}}{k!}\right) dx$$

$$= \sum_{k=0}^{[n/2]} \frac{(-[n/2])_k}{(b+(-1)^n/2)_k\, k!} \left(\int_{-\infty}^{\infty} e^{-isx} e^{\frac{-1}{2x^2}} x^{-2a+n-2k}\,dx\right).$$
(41)

Again, the following definite integral should be evaluated

$$R_{n,k}(s;a) = \int_{-\infty}^{\infty} e^{-isx} e^{\frac{-1}{2x^2}} x^{-2a+n-2k}\,dx. \quad (42)$$

To do this, we can use the same method as we applied to $I_{n,k}(s;\alpha,\beta)$, i.e.

$$R_{2m,k}(s;a) = \int_{-\infty}^{\infty}\left(\sum_{j=0}^{\infty} \frac{(-isx)^j}{j!}\right) e^{\frac{-1}{2x^2}} x^{-2a+2m-2k}\,dx = \sum_{j=0}^{\infty} \frac{(-1)^j i^j s^j}{j!}\left(\int_{-\infty}^{\infty} e^{\frac{-1}{2x^2}} x^{j-2a+2m-2k}\,dx\right)$$

$$= \sum_{r=0}^{\infty} \frac{(-1)^r s^{2r}}{(2r)!}\left(2\int_0^{\infty} x^{2r-2a+2m-2k} e^{\frac{-1}{2x^2}}\,dx\right) = \sum_{r=0}^{\infty} \frac{(-1)^r s^{2r}}{(2r)!} 2^{-r+a-m+k-\frac{1}{2}} \Gamma(-r+a-m+k-\frac{1}{2})$$

$$= 2^{a-m+k-\frac{1}{2}} \Gamma(a-m+k-\frac{1}{2})\, _0F_2\left(\begin{matrix}-\\ \frac{1}{2},\ \frac{3}{2}-a+m-k\end{matrix}\bigg|\frac{s^2}{8}\right),$$
(43)

as well as

$$R_{2m+1,k}(s;a) = \sum_{j=0}^{\infty} \frac{(-1)^j i^j s^j}{j!} \left( \int_{-\infty}^{\infty} e^{\frac{-1}{2x^2}} x^{j-2a+2m+1-2k} dx \right) = (-is) \sum_{r=0}^{\infty} \frac{(-1)^r s^{2r}}{(2r+1)!} \left( 2 \int_{0}^{\infty} x^{2r-2a+2m-2k+2} e^{\frac{-1}{2x^2}} dx \right)$$

$$= (-is) 2^{a-m+k-\frac{3}{2}} \Gamma(a-m+k-\frac{3}{2}) \, _0F_2\left( \begin{array}{c} - \\ \frac{3}{2}, \frac{5}{2}-a+m-k \end{array} \bigg| \frac{s^2}{8} \right).$$

(44)

Consequently we have

$$R_{n,k}(s;a) = 2^{a+k-\frac{1}{2}-[\frac{n+1}{2}]} \Gamma\left(a+k-\frac{1}{2}-[\frac{n+1}{2}]\right)(-is)^{\frac{1-(-1)^n}{2}} \, _0F_2\left( \begin{array}{c} - \\ 1-\frac{(-1)^n}{2}, -a-k+\frac{3}{2}+[\frac{n+1}{2}] \end{array} \bigg| \frac{s^2}{8} \right),$$

(45)

and therefore

$$\mathbf{F}(u(x)) = \Gamma(a-\frac{1}{2}-[\frac{n+1}{2}]) 2^{a-\frac{1}{2}-[\frac{n+1}{2}]} (-is)^{\frac{1-(-1)^n}{2}} \times$$

$$\sum_{k=0}^{[n/2]} \frac{(-[\frac{n}{2}])_k (a-\frac{1}{2}-[\frac{n+1}{2}])_k}{(b+(-1)^n/2)_k} \frac{2^k}{k!} \, _0F_2\left( \begin{array}{c} - \\ 1-\frac{(-1)^n}{2}, -a-k+\frac{3}{2}+[\frac{n+1}{2}] \end{array} \bigg| \frac{s^2}{8} \right).$$

(46)

Now if for simplicity in (46) we define the symmetric function

$$B_n(x;q_1,q_2) = x^{\frac{1-(-1)^n}{2}} \sum_{k=0}^{[n/2]} \frac{(-[\frac{n}{2}])_k (q_1-\frac{1}{2}-[\frac{n+1}{2}])_k}{(q_2+(-1)^n/2)_k} \frac{2^k}{k!} \, _0F_2\left( \begin{array}{c} - \\ 1-\frac{(-1)^n}{2}, -q_1-k+\frac{3}{2}+[\frac{n+1}{2}] \end{array} \bigg| \frac{x^2}{8} \right),$$

(47)

then by referring to definitions (40) and applying Parseval identity we get

$$2\pi \int_{-\infty}^{\infty} x^{-2(a+c)} e^{\frac{-1}{x^2}} B_n^{(b)}(x) B_m^{(d)}(x) dx$$

$$= i^{\frac{(-1)^n-(-1)^m}{2}} \frac{\Gamma(a-\frac{1}{2}-[\frac{n+1}{2}])\Gamma(c-\frac{1}{2}-[\frac{m+1}{2}])}{2^{-a+\frac{1}{2}+[\frac{n+1}{2}]} 2^{-c+\frac{1}{2}+[\frac{m+1}{2}]}} \int_{-\infty}^{\infty} B_n(s;a,b) B_m(s;c,d) ds.$$

(48)

It is sufficient in (48) to assume that $b = d = a+c$ and then refer to the finite orthogonality relation (17) to reach,

**Theorem 2.** *The special function $B_n(x;q_1,q_2)$ defined in (47) satisfies the orthogonality relation*

$$\frac{1}{2\pi}\int_{-\infty}^{\infty} B_n(x;a,b)B_m(x;b-a,b)\,dx = 2^{-b+1+2[\frac{n+1}{2}]}\prod_{j=1}^{n}\frac{2(-1)^j(j-b)+2b}{(2j-2b+1)(2j-2b-1)}$$
$$\times \frac{\Gamma(b-1/2)}{\Gamma(a-\frac{1}{2}-[\frac{n+1}{2}])\Gamma(b-a-\frac{1}{2}-[\frac{n+1}{2}])}\delta_{n,m}, \quad (49)$$

for $m,n = 0,1,\ldots,N = \max\{m,n\} \leq b - \frac{1}{2}$, $(-1)^{2b} = 1$ and $\frac{1}{2} < a < b - \frac{1}{2}$.